\documentclass{article}
\usepackage{amsmath,graphicx,amssymb,amsfonts,amsthm}

\usepackage{epstopdf}
\title{Blind calibration for compressed sensing by convex optimization}
\author{R. Gribonval, G. Chardon and L. Daudet}
\def\y {\mathbf{y}}
\def\Y {\mathbf{Y}}
\def\x {\mathbf{x}}
\def\X {\mathbf{X}}
\def\D {\mathbf{D}}
\def\M {\mathbf{M}}

\def\diag {\mathrm{diag}}
\def\R {\mathbb{R}}
\def\C {\mathbb{C}}
\def\n {\mathbf{n}}
\def\N {\mathbf{N}}
\def\estM {\widehat{\M}}
\def\tM {\widetilde{\M}}
\def\estD {\widehat{\D}}
\def\estX {\widehat{\X}}
\def\iD {\mathbf{\Delta}}
\def\estiD {\widehat{\iD}}
\def\Tr {\mathrm{Tr}}

\begin{document}
\maketitle

\begin{abstract}
We consider the problem of calibrating a compressed sensing measurement system under the assumption that the decalibration consists in unknown gains on each measure. We focus on {\em blind} calibration, using measures performed on a few unknown (but sparse) signals. A naive formulation of this blind calibration problem, using $\ell_{1}$ minimization, is reminiscent of blind source separation and dictionary learning, which are known to be highly non-convex and riddled with local minima. In the considered context, we show that in fact this formulation can be exactly expressed as a convex optimization problem, and can be solved using off-the-shelf algorithms.  Numerical simulations demonstrate the effectiveness of the approach even for highly uncalibrated measures, when a sufficient number of (unknown, but sparse) calibrating signals is provided. We observe that the success/failure of the approach seems to obey sharp phase transitions.
\end{abstract}

%\begin{keywords}
%compressed sensing; calibration; dictionary learning; blind signal separation; s%parse recovery.
%\end{keywords}

\section{Introduction}

Linear inverse problems are ubiquitous in signal and image processing, where they are used to estimate an unknown signal $\x_{0} \in \R^{N}$ or $\C^{N}$ from noisy linear measurements:
\[
\y := \M \x_{0} + \n \in \R^{m}\ \mbox{or} \ \C^{m}.
\]
When $m<N$ this well-known under-determined problem admits infinitely many solutions, but if $\x_{0}$ is sparse enough it can be estimated accurately using sparse regularization. Among many other techniques, $\ell_{1}$ regularization has recently become quite popular and consists in solving, e.g.,
\[
\hat\x_{1} := \arg\min_{\x} \|\x\|_{1}\ \mbox{s.t.}\ \|\y-\M\x\|_{2} \leq \epsilon
\]
for a well chosen $\epsilon$. 
When the measurement matrix $\M$ is perfectly known, this approach is known to perform well, and a well-established body of work characterizes its performance guarantees for the recovery of vectors $\x_{0}$ that are sufficiently well approximated by highly sparse vectors (see, e.g.~\cite{fuchs04:_spars_repres_arbit_redun_bases,Candes:2008aa,Baraniuk:2010aa}). 
Such sparse linear regression problems occur in many practical scenarii where the measurement matrix $\M$ is either dictated by the physics of the measurement system, or designed to have favorable properties with respect to the recovery of sparse vectors: this is the now famous compressed sensing (CS) scenario (see, e.g.~\cite{MR2241189,stablesigrecovery-CandesRombergTao-2006,Baraniuk:2007aa}), where $\M$ {\em voluntarily} reduces dimension, exploiting the sparsity of $\x_{0}$ to capture it with fewer measurements than Nyquist sampling would require.

\subsection{The decalibration issue}
In practical situations, the true measurement system $\M$ is not perfectly known: it may only have been modeled; or it may have been measured through a calibration process, but the physical conditions of the system (such as temperature) may have drifted since this calibration. 

Exploiting sparse regularization 
%\[
%\hat\x_{1} := \arg\min_{\x} \|\x\|_{1}\ \mbox{s.t.}\ \|\y-\estM\x\|_{2} \leq \epsilon
%\]
with an inaccurate estimate $\estM$ of the true measurement system $\M$ is likely to hurt the reconstruction performance~\cite{5419044}. It is believed to be one of the reasons limiting the observed performance of several compressed sensing devices.
To address this problem, the most standard existing approaches are: 
\begin{enumerate}
\item {\bf To ignore the problem.}
\item {\bf To consider de-calibration as noise~\cite{5419044}:} $\y \approx \estM \x_{0} + (\M-\estM)\x_{0}$. This leads to solving
\[
\hat\x_{1} := \arg\min_{\x} \|\x\|_{1}\ \mbox{s.t.}\ \|\y-\estM\x\|_{2} \leq \epsilon+\eta
\]
with $\eta$ an estimate of the magnitude of this added noise.
\item {\bf Supervised calibration:} using known training signals $\x_{1}, \ldots \x_{L}$ and the corresponding observations $\y_{\ell} = \M \x_{\ell} + \n_{\ell}$. Gathering all data in matrices, this takes the form $\Y = \M \X + \N$. The matrix $\M$ is re-estimated, e.g., as
\[
\estM := \arg\min_{\tM} \|\Y-\tM\X\|_{F}^{2}.
\]
\end{enumerate}

\subsection{Constraints on calibration}
It is sometimes useful to constrain the estimated calibration matrix $\estM$ to belong to some family $\mathcal{M}$ of matrices. 
For example, it is sometimes known/assumed~\cite{Pfander:2008aa} that the unknown $\M$ is sparse in a given dictionary of measurement matrices $\{\M_{k}\}$: $\M \approx \sum_{k} \alpha_{k} \M_{k}$, $\|\mathbf{\alpha}\|_{0}$ small. Supervised calibration can then be performed by solving, e.g., the convex relaxation
\[
\min \|\mathbf{\alpha}\|_{1} \ \mbox{s.t.}\ \|\Y - \sum_{k} \alpha_{k} \M_{k} \X\|_{F} \leq \epsilon.
\]

In this paper, we concentrate on a different scenario where the measurement matrix is almost known, up to an unknown gain on each measure: that is to say, $\M = \D_{0} \M_{0}$ where $\M_{0}$ is a perfectly known measurement matrix, and $\D_{0}$ is an unknown diagonal matrix which $i$-th entry is a (real or complex) gain $d_{i}$ applied to the $i$-th measure of $\x_{0}$ associated to the $i$-th line of $\M_{0}$. This leads to 
\[
\mathcal{M} := \{\M = \D \M_{0}, \D = \diag(d_{i}), d_{i} \neq 0\ \forall i\}.
\]
Several practical scenarii can be associated with this assumption, for instance in the case of a microphone array where the frequency response of each microphone needs to be individually calibrated~\cite{Mignot:2011aa}. At each frequency, the calibration problem amounts to choosing unknown gains.

\subsection{Unsupervised / blind sparse calibration}
In this paper, we are interested in a {\em blind} calibration problem, where no {\em known} training signal can be used. While some training signals $\x_{\ell}$ have given rise to observed measures $\y_{\ell}$, the training signals themselves are not known. Since one could not hope to calibrate the system without some form of knowledge of $\x_{\ell}$, and the final scenario for using $\estM$ will be sparse reconstruction, we assume that the unknown signals $\x_{\ell}$ are sparse and somehow statistically diverse. 
In matrix form, we have
\(
\Y = \D_{0} \M_{0} \X_{0}
\)
where $\Y$ is the $m \times N$ (known) observation matrix, $\D_{0}$ is the square and diagonal matrix of size $m$ for the (unknown) calibration coefficients, $\M_{0}$ is the $m \times N$ (known) idealized measurement matrix, and $\X_{0}$ is the (unknown) $N \times L$ set of cali\-bra\-tion signals upon which the only knowledge we have is that they are $k$-sparse.

The objective is to derive a calibration technique that exploits the knowledge of $\Y$ and $\M_{0}$, as well as the sparsity of the unknown $\X_{0}$, to obtain estimates $\estD$ and $\estX$, and finally $\estM = \estD \M_{0}$. It seems natural to consider the following optimization problem, where the $\ell_{1}$ objective function is intended to promote the sparsity of $\X$:
\begin{equation}
\label{eq:BlindCalibration}
\min_{\D,\X} \|X\|_{1}\ \mbox{s.t.} \Y = \D \M_{0}\X.
\end{equation}
However, this would naively lead to two major issues: 
\begin{enumerate}
\item Without further constraint on $\D$ and/or $\X$, one can scale $\D$ to infinity while letting $\X$ go to zero, leading to a trivial but uninteresting solution;\\ This is usually solved through a normalization constraint $\D \in \mathcal{D}$ where the set $\mathcal{D}$ is bounded;
\item Even with appropriate normalizing constraints on $\D$, this seems at first a non-convex problem, because of the bilinear nature of the term $\D \M_{0}\X$, which is linear / convex separately in $\D$ and $\X$ but not jointly. 
\end{enumerate}

\subsection{Relation to previous work}
The above issues are well known since they are encountered in {\em blind source separation}~\cite{Zibulevsky:2001aa} and {\em dictionary learning}~\cite{2003_NeuralComp_DictLearning_KreutzEtAl,2006_IEEE_TSAP_KSVD_AharonEtAl}. In a way, the unsupervised calibration problem can be seen as a simplified instance of dictionary learning, where the general problem has been expressed, e.g, as~\cite{Yaghoobi:2008aa,4787130,Gribonval:2010aa,Mairal:2010aa}
\[
\min_{\mathbf{\Phi},\X} \|\X\|_{1}\ \mbox{s.t.} \Y = \mathbf{\Phi} \X.
\]
where the minimum is over $\mathbf{\Phi} \in \mathcal{C}$ for a well-chosen bounded set $\mathcal{C}$ (typically, the set of matrices with unit norm columns, or {\em oblique manifold}), but the dictionary $\mathbf{\Phi}$ is not restricted to be diagonal.

%\subsection{Link with the  problem}
%Maybe this is where we could discuss the differences and similarities. 
  The present issue is also related to the {\em basis mismatch} problem~\cite{chi2009sensitivity}, where the signals are truly sparse in a basis that is slightly different from the chosen representation basis (for instance with a parametric dictionary using a grid of quantized parameters); or to the perturbed CS problem~\cite{herman2010general, zhu2010sparsity}, where multiplicative noise limits the effectiveness of CS. In both cases, the ``noise'' on the measurement matrix is different for every vector, whereas in the  decalibration case considered here, this (unknown) basis is shared by all vectors. 

\subsection{Contributions}
The main contribution of this paper is (Section~\ref{sec:Approach}) to show that, under proper parametrization and normalization, the considered unsupervised calibration problem~\eqref{eq:BlindCalibration} can actually be {\em exactly} expressed as a convex optimization problem. 
Numerically solving this problem is straightforward using off-the-shelf algorithms. In numerical experiments (Section~\ref{sec:Numerics}), we demonstrate the effectiveness of this approach even for highly uncalibrated measurements, whenever a sufficient number of (unknown but sparse) calibrating signals is provided. Remarkably, we observe that the success/failure of this problem obeys some sharp phase transitions, generalizing to the uncalibrated case the phase transitions studied by Donoho and Tanner~\cite{Donoho:2009aa}.

\begin{figure*}[htbp] %  figure placement: here, top, bottom, or page
   \centering
     \includegraphics[width=14 cm]{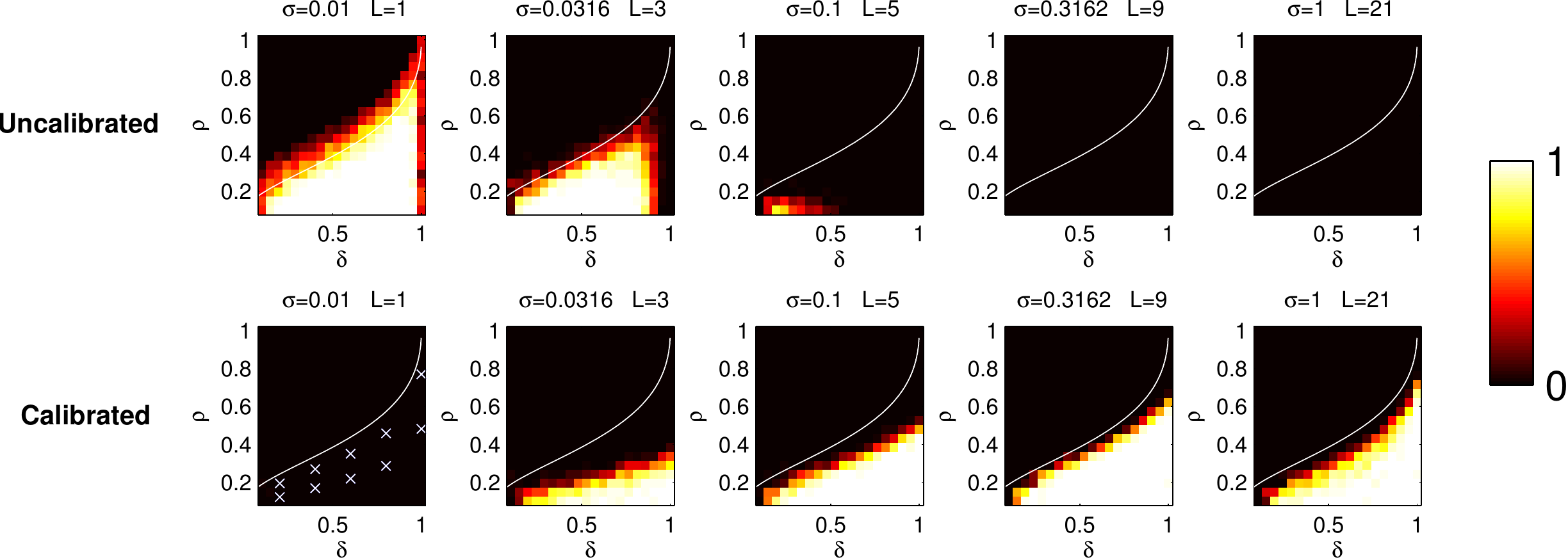}
   \caption{Empirical phase-transitions with no calibration (top) and proposed calibration (bottom). Overwhelming success (white) or failure (black) is displayed as a function of dimension parameters $\delta := m/N$ (abscissa) and $\rho := k/m$ (ordinate). From left to right: increasing values of the number $L$ of training signals and decalibration amplitude $\sigma$. The thin line is the asymptotic Donoho-Tanner phase transition curve for known Gaussian measurement matrix. Crosses on the bottom-left figure indicate choices of  ( $\delta$, $\rho$ ) for  the experiments on Fig.~\ref{fig:NScurve}.}
   \label{fig:DTcurve}
\end{figure*}

\section{Proposed approach}
\label{sec:Approach}
As noticed above, the naive formulation of the problem is non-convex, even with a (convex) normalization constraint. 
\subsection{Convex formulation}
To provide a convex formulation we propose to reparameterize the problem. Using the simple assumption that the unknown gains $d_{i}$ are nonzero, we can write the constraint $\Y = \D \M_{0} \X$ as $\iD \Y = \M_{0} \X$, where $\iD = \D^{-1} = \diag(\delta_{i})$. 
\subsection{Normalization constraint}
With this new parametrization, the constraint is convex in the pair $(\iD,\X)$. However, it is trivially satisfied for the pair $(0,0)$. To avoid this trivial solution we introduce a convex normalization constraint $\Tr(\iD) = \sum_{i} \delta_{i} = m$. Therefore, the constraint set is $\mathcal{D} := \{\D = \diag(d_{i}),\ \sum_{i} d_{i}^{-1} = m\}$. Note that many other alternatives could be considered such as $\delta_{1} = 1$. 
\subsection{Proposed blind calibration approach}
We end up proposing the following unsupervised calibration approach: given the collection of observed measures $\Y$ and the model of the measurement system $\M_{0}$, we estimate the inverse of the calibration gains $\estiD$ and the training signals $\estX$ as:
\begin{equation}
\label{eq:Calibrated}
(\estX_{\textrm{cal}},\estiD) := \arg\min_{\X,\iD} \|\X\|_{1} \ \mbox{s.t.}\ \iD\Y = \M_{0}\X,\ \Tr(\iD)=m.
\end{equation}
This is a convex problem.

\section{Experimental results}
\label{sec:Numerics}

\subsection{Considered techniques}
We compare two techniques. The proposed approach~\eqref{eq:Calibrated}, and an approach which ignores decalibration
\[
\estX_{\textrm{uncal}} := \arg\min_{\X} \|\X\|_{1} \ \mbox{s.t.}\ \Y = \M_{0}\X.
\]
Both convex problems are solved using the \texttt{cvx} toolbox~\cite{cvx}.

\subsection{Data generation}

The data was generated by drawing i.i.d. random $k$-sparse vectors $\x_{\ell}$ with a $k$-sparse support chosen uniformly at random and i.i.d. Gaussian nonzero entries. The idealized measurement matrix $\M_{0}$ was drawn from the Gaussian ensemble of size $m \times N$ with $N = 100$. Real positive decalibration coefficients were generated using   
%$\D_{0} = \Id + \sigma \cdot \textrm{diag}(\mathcal{N}(0,1))$.
%should replace with
$\D_{0} = \diag(\exp(\mathcal{N}(0,\sigma^{2})))$, where $\sigma$ is the  parameter governing the amplitude of decalibration. 
Hence, a given value of $\sigma$ leads to a decalibration offset in decibels with zero mean and standard deviation $\pm 20 \sigma / \ln(10) \approx  \pm 8.7  \sigma$ dB. 
%This exponential calibration is on par with the "physics" (off by a certain number of dB). 
The considered measures were $\Y = \D_{0}\M_{0}\X$, where $\X$ was a concatenation of $L$ i.i.d. $k$-sparse signals.
Experiments were conducted for various configurations of $\delta := m/N$, $\rho := k/m$, $L$ and $\sigma$ to measure how accurately the signals $\X$ were recovered with the different approaches. Recovery was considered successful if the normalized cross-correlation between original and estimated signals is above 99.5 \% (note that, due to the global scaling invariance in $\iD$ and $\X$ mentioned in the introduction, a distance measure between $\X$ and $\estX$ would not be a relevant indication of success).
All the empirical phase transition diagrams described below report the number of successful estimations, when the experiment was repeated over 50 random draws. 

%\subsection{Performance measure (skip if we are lacking space ?)}
%Since $\Tr(\D_{0})$ is a priori unknown and might differ from $m$ (i.e., all calibrations coefficients $d_i$ average to 1), the solution to the considered calibration problem may suffer from a global scaling. As a result, instead of measuring the calibration performance by $\|\estX-\X_{0}\|_{F}$, we rather consider a performance measure based on correlation:
%\[
%\textrm{Perf}(\estX,\X_{0}) := \frac{\langle \estX,\X_{0} \rangle_{F}} {\|\estX\|_{F} \cdot \|\X_{0}\|_{F}}
%\]
%where by definition $\langle \estX,\X_{0} \rangle_{F} := \Tr(\estX^{T}\X_{0})$ is the inner product associated to the Froebenius norm $\|\X\|_{F}^{2} = \langle \X,\X \rangle_{F} = \Tr \X^{T}\X = \Tr \X\X^{T}$.

\subsection{Can CS still succeed under calibration errors ?}
The first study dealt with estimating how much classical CS was robust to calibration errors. To do this, we ran numerical experiments and plotted empirical Donoho-Tanner phase transition~\cite{Donoho:2009aa} (in the framework presented above). 
Different de-calibration levels $\sigma$ and number $L$ of calibration signals have been tested. 
%These plots show the percentage of success for the reconstruction, for different oversampling ratio 
%$\delta := m/N$, and sparsity ratio $\rho := k/m$.

The results can be seen on the top line of Fig.\ref{fig:DTcurve}:  for small values of de-calibration ($\sigma = .01$, $\sigma = .0316$), the transition curve is barely modified : CS seems indeed robust to small calibration errors. However, as the de-calibration increases ($\sigma = .1$, i.e. the decalibration error was of the order $\pm .9$ dB), then the region where CS succeeds drastically shrinks, and eventually disappears at even larger de-calibration ($\sigma = .316$ and above).  

The second line of fig. \ref{fig:DTcurve} shows the result of CS under the blind calibration procedure introduced above. For a large number of noise levels where the un-calibrated experiment failed, it now succeeds even with a relatively small number $L$ of training samples (for instance, with  $L = 21$ training samples, it still succeeds at $\sigma = 1$, i.e., a decalibration error of the order $\pm 8.7$ dB).    
Interestingly, for very small values of de-calibration ($\sigma = .01$ and below) and too few training samples, it is better not to perform the blind calibration, as it introduces too many degrees of freedom. 
\begin{figure*}[htbp] %  figure placement: here, top, bottom, or page
   \centering
\includegraphics[width=14 cm]{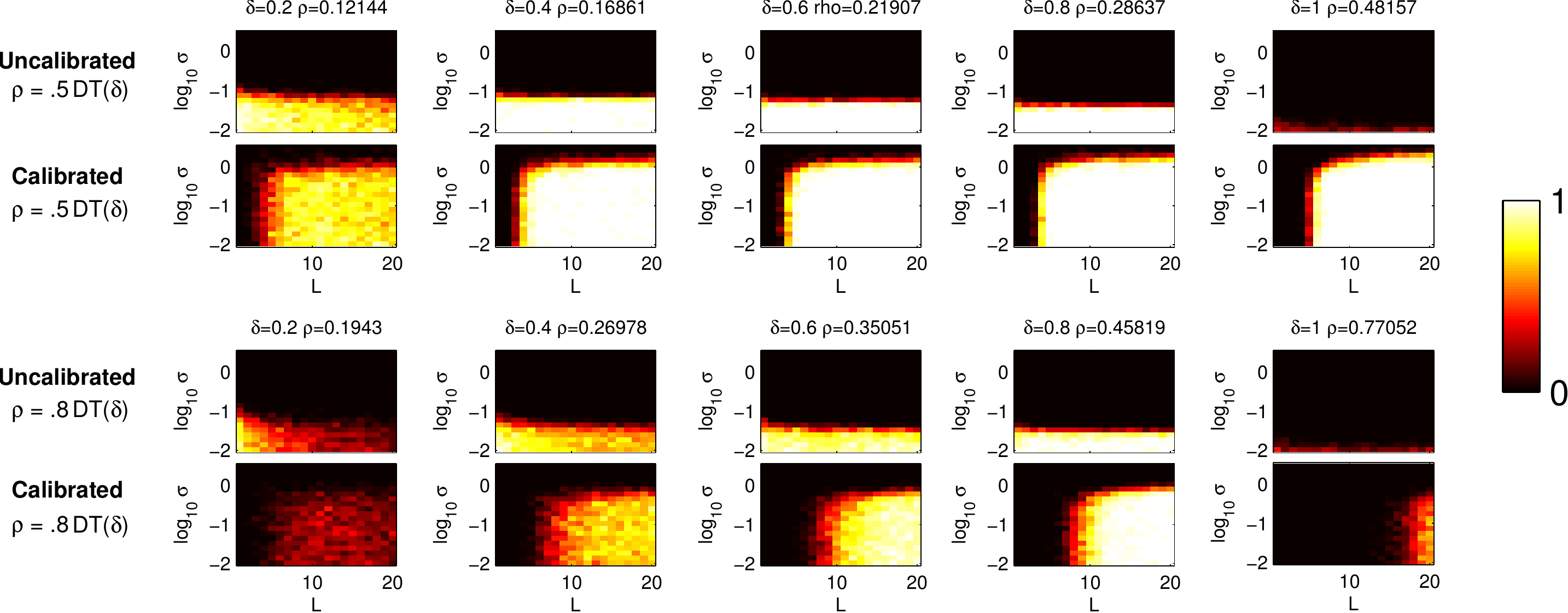}%matlab/NS_article2.pdf} 

   \caption{Experiments at fixed ($\delta$,$\rho$), as a function of the number $L$ of training signals and the decalibration amplitude $\sigma$. ``Uncalibrated'' and ``calibrated'' refer to the standard and proposed  algorithms, respectively. White indicates that the algorithm succeeds in recovering the solution, and black indicates a failure.    }
   %[we should probably just keep the first 2 lines corresponding to rho = 0.5 DT(delta), the other ones are for .8 DT]
   \label{fig:NScurve}
\end{figure*}
%\marginpar{We should put on one of the DT images the location of $\delta,\rho$ considered in Fig. 2? On figures, replace Ntrain by L}
 
\subsection{Choosing the number of training samples}
In this section we wish to determine how many training samples $L$ are needed as a function of the ``noise'' level $\sigma$. We picked different values of $(\delta, \rho)$ where 
ideal CS works, and plotted on fig.\ref{fig:NScurve} the rate of success as a function of de-calibration level $\sigma$ and number of training samples $L$, both in the un-calibrated and calibrated cases. Remarkably, these obey sharp corner-like phase transition. Calibration allows CS to work for $\sigma$ increased by an order of magnitude (compared to the un-calibrated case), provided $L$ is larger than some threshold. These threshold values seem to have some universality, as numerical tests do not indicate any clear dependency on the ambient dimension $N$ (results not shown here), and sharper transitions are observed for larger values of $N$.

\section{Conclusion}
% {\bf a few words of conclusion?}
In this study, we have experimentally observed that even mild decalibration on the measurement matrix can lead to a spectacular failure of standard CS recovery algorithms. This problem is highly relevant for engineering applications, as in many cases the observation matrix is not exactly known. 
The proposed blind calibration technique, formulated here as a convex optimization problem, is shown to offer a significantly improved robustness to decalibration: when a sufficient number of (unknown but sparse) calibration signals are provided, the algorithm can now succeed even for much larger calibration errors. Remarkably, the success / failure of this new calibrated CS exhibits sharp corner-like phase transitions. 
Being able to compute theoretically the asymptotic boundaries of these new transition diagrams is a promising extension of this work, together with the extension to more general decalibration cases. More immediate targets will be the design of efficient blind calibration algorithms to scale the approach and test in on real calibration problems arising in compressed sensing of acoustic fields~\cite{Mignot:2011aa}.

%For $\sigma = 0.1$ we tried $L=1:10$ and observed that for $L= 5$ the achieved SNR with $\hat\X_{2}$ exceeds $100$ dB.

%Changed $N$ in $L$ ($N$ is the dimension of $x$). 

%\section{}
%\subsection{}
%\ninept
%\bibliographystyle{plain}
\bibliographystyle{IEEEbib}
\bibliography{../../../../BiblioDatabaseASCII,../../../../Gribonval}

\end{document}